\documentclass[11pt]{amsart}
\usepackage{bm}
\usepackage{amssymb}
\usepackage{graphicx}	
\usepackage{xcolor}
\usepackage{amsmath}
\usepackage{enumerate}
\usepackage{blindtext}
\usepackage{multirow}
\usepackage{mathrsfs}

\newtheorem{thm}{Theorem}[section]

\newtheorem{cor}[thm]{Corollary}

\newtheorem{Def}[thm]{Definition}
\newtheorem{prop}[thm]{Proposition}
\newtheorem{rem}[thm]{Remark}
\newtheorem{ex}[thm]{Example}

\newcommand{\bdfn}{\begin{Def} \rm}
\newcommand{\edfn}{\end{Def}}

\newcommand{\tfae}{the following are equivalent}

\newcommand{\ra}{\rightarrow}

\newcommand{\Lglra}{\Longleftrightarrow}
\newcommand{\es}{\emptyset}
\newcommand{\vp}{\varphi}
\newcommand{\ci}{\subseteq}

\newcommand{\al}{\alpha}

\newcommand{\de}{\delta}
\newcommand{\e}{\varepsilon}

\newcommand{\la}{\lambda}
\newcommand{\si}{\sigma}
\newcommand{\otm}{\otimes}
\newcommand{\ga}{\gamma}
\newcommand{\Ti}{\Tilde}

\newcommand{\mb}{\mathbb}
\newcommand{\mc}{\mathcal}
\newcommand{\tr}{\textrm}
\newcommand{\mf}{\mathfrak}

\newcommand{\mr}{\mathscr}

\newcommand{\Si}{\Sigma}

\newcommand{\iy}{\infty}
\newcommand{\om}{\omega}
\newcommand{\Om}{\Omega}

\newcommand{\bo}{\bigoplus}

\newcommand{\trC}{\textrm{Cent}}
\newcommand{\trr}{\textrm{rad}}

\newcommand{\ercp}{\bf{\emph{rcp}}}
\newcommand{\eSACP}{\bf{\emph{SACP}}}

\newcommand{\beqa}{\begin{eqnarray*}}
\newcommand{\eeqa}{\end{eqnarray*}}
\newcounter{cnt1}
\newcounter{cnt2}
\newcounter{cnt3}
\newcounter{cnt4}
\newcommand{\blr}{\begin{list}{$($\roman{cnt1}$)$} {\usecounter{cnt1}
\setlength{\topsep}{0pt} \setlength{\itemsep}{0pt}}}
\newcommand{\blR}{\begin{list}{\Roman{cnt4}.\ } {\usecounter{cnt4}
\setlength{\topsep}{0pt} \setlength{\itemsep}{0pt}}}
\newcommand{\bla}{\begin{list}{$(\alph{cnt2})$} {\usecounter{cnt2}
\setlength{\topsep}{0pt} \setlength{\itemsep}{0pt}}}
\newcommand{\bln}{\begin{list}{$($\arabic{cnt3}$)$} {\usecounter{cnt3}
\setlength{\topsep}{0pt} \setlength{\itemsep}{0pt}}}
\newcommand{\el}{\end{list}}

\DeclareMathOperator{\spn}{span}
\newcommand{\GC}{{\mbox{\rm (GC)}}}

\begin{document}
\title[Local constraints and Best Simultaneous Approximation]{Transfer of Approximation properties under Local Constraints and Best Simultaneous Approximation on Sums}

\author[Das]{Syamantak Das}

\newcommand{\Addresses}{{
  \bigskip
  \footnotesize
  Syamantak Das, \textsc{National Institute of Technology Sikkim, India }\par\nopagebreak
  \textit{E-mail address}, Syamantak Das: \texttt{syamantakdas@nitsikkim.ac.in}
 
}}

\subjclass[2000]{Primary 41A28, 41A65 Secondary 46B20 41A50 \hfill \textbf{\today} }
\keywords{constrained subspace, central subspace, simultaneous approximation, generalized center}
\begin{abstract}

It is folklore that the sum of two $M$-ideals (semi $M$-ideals) is also an $M$-ideal (a semi $M$-ideal). Numerous authors have attempted to investigate such properties of subspaces. This article explores two important aspects of approximation theory within Banach spaces and how these properties remain intact when considering the sum of two subsets. Recall the notion of $(GC)$ introduced by Vesel\'y that encloses two aforementioned properties. When the sum of two subspaces is closed, we discuss various properties of the sum if one of the subspaces has these properties. Counterexamples are produced that establish nonaffirmativeness for the properties $(GC)$ and the central subspace.
We answer a problem raised by the author in [{\em Best constrained approximation in Banach spaces}, Numer. Funct. Anal. Optim. {\bf 36}(2) (2015),  248--255]. We extend our observations related to the best simultaneous approximations to the properties $(P_1)$ and $\mr{F}$-SACP.

\end{abstract}
 \maketitle
\section{Introduction}

In this paper we focus on the transfer of several approximation properties under local constraints and then direct our attention to the sum of two subsets having these properties.
Our first main result (Section~2) settles a question posed by Rao by showing that under local constraints the $M$-ideal property descends from $Z_1\subset X$ to $Z_2=Z_1\cap Y$ in $Y$.
Abstracting from its proof a unified transfer mechanism, we obtain stability results for central subspaces and $(\mathrm{GC})$, with functorial extensions to polyhedral direct sums, K\"othe--Bochner spaces, and spaces obtained by injective tensor product of spaces.
In a complementary direction (Section~3), we develop a gauge-based framework for best simultaneous approximation and establish transfer theorems for $\mr{F}$-rcp, $\mr{F}$-property-$(P_1)$, and ($\mr{F}$-SACP) on closed sums.

\subsection{Notations and Preliminaries}

$X$ denotes a real Banach space. By a subspace $Y$ of $X$ we indicate a closed linear subspace of $X$. $B_X$ denotes the closed unit ball of $X$. For $x\in X$ and $r>0$, $B_X(x,r)$ denotes the closed ball in $X$ centered at $x$ and radius $r$. When there is no chance of confusion, we simply write it as $B(x,r)$. We consider an $x\in X$ to be canonically embedded in $X^{**}$. 

A subspace $Y$ of $X$ is said to be a {\it constrained subspace} (or {\it 1-complemented}) if $Y$ is the range of a norm-1 projection on $X$. Constrained subspaces of higher dimension are rare in Banach spaces. In certain instances, obtaining a local version, specifically for all $x\notin Y$, the existence of a norm-$1$ projection from $\spn\{x,Y\}$ onto $Y$ is also rare. For example, $c_0$ in $\ell_\iy$ is one such case. A subspace $Y$ of $X$ such that for all $x\notin Y$, there exists a norm-1 projection $P:\spn\{x,Y\}\ra Y$ is said to be an {\it almost constrained subspace} of $X$. In Chapter~5 in \cite{JL1} Lindenstrauss studies these subspaces.

In \cite{L} Lindenstrauss gives an example of a Banach space $X$ and a subspace $Y$ of $X$ such that dim$(X/Y)=2$ and for all $x\notin Y$ there is a norm-1 projection from $\spn\{x,Y\}$ onto $Y$ but $Y$ is not a constrained subspace of $X$. 

Several weakenings of a constrained subspace, such as almost constrained subspace, $\GC$, central subspace, locally constrained subspace have been studied in the literature. Ves\'ely in \cite{LV} introduces the notion of generalized centers in Banach spaces. Recall the following result from \cite{LV}. For any unexplained terminology, we refer to the subsequent part of this section and also the paper \cite{LV}. 
\begin{thm}\label{T4}
For a finite subset $\{x_1,\cdots,x_N\}$ of $X$, \tfae.
\bla
\item If for $r_1,\cdots,r_N>0$, $\cap_{i=1}^NB_{X^{**}}(x_i,r_i)\neq\es$, then $\cap_{i=1}^NB_X(x_i,r_i)\neq\es$.
\item $\{x_1,\cdots,x_N\}$ admits weighted Chebyshev centers for all weights $\rho_1,\cdots,\rho_N>0$.
\item $\{x_1,\cdots,x_N\}$ admits $f$-centers for every continuous monotone coercive function $f$ on $\mb{R}^N_+$.
\el
\end{thm}

$X\in \GC$ if for any finite subset $\{x_1,\cdots,x_N\}$ of $X$, it satisfies any one of the equivalent conditions of Theorem~\ref{T4}. We note that the condition $(a)$ of Theorem~\ref{T4} is in some sense a strengthening of the principle of local reflexivity in Banach spaces.

Motivated by the notion of $\GC$, Bandyopadhyay and Rao introduce the notion of a central subspace in \cite{PB}.

\bdfn
A subspace $Y$ of $X$ is said to be a {\it central subspace} of $X$ if for any finite family of closed balls with centers in $Y$ having a nonempty intersection in $X$ also intersect in $Y$.
\edfn

Bandyopadhyay and Dutta further generalized the notion of a central subspace in \cite{SD} to an {\it $\mc{A}$-$C$-subspace}, for a subfamily $\mc{A}$ of subsets of $X$.

\bdfn\cite{SD}
Let $Y$ be a subspace of $X$  and $\mc{A}$ be a family of subsets of $Y$.  $Y$ is said to be an {\it $\mc{A}$-$C$-subspace of $X$} if for  $x\in X$ and $A\in\mc{A}$, there exists $y\in Y$ such that $\|y-a\|\leq\|x-a\|$ for all $a\in A$.

If $\mc{A}$ is a family of subsets of $X$, then it is referred to that  $X$ has {\it $\mc{A}$-IP}, if $X$ is an $\mc{A}$-$C$-subspace of $X^{**}$.
\edfn

Rao in \cite{Rao} introduces the notion of locally constrained subspaces.
\bdfn\label{D51}
Let $Z_1,Y$ be subspaces of $X$ and $Z_2=Z_1\cap Y$.  $Z_2$ is said to be {\it locally constrained with respect to $Z_1$ via $Y$}, provided  for all $z\in Z_1$, there are linear projections $P:\spn\{z,Z_2\}\longrightarrow Z_2$ and $Q:\spn\{z,Y\}\longrightarrow Y$ of norm one, such that $Pz=Qz$.
\edfn

In \cite{Rao}, it is observed, in particular, that if $Y$ is a constrained subspace of $X$, then $Z_2$ is locally constrained with respect to $Z_1$ via $Y$. Moreover, if $Y$ and $X$ are as stated above in the example by Lindenstrauss, then $Y$ is locally constrained with respect to $X$ via $Y^{**}$, is observed in \cite{Rao}.

The study on the existence of generalized centers in $X$ (and more generally on a closed subset $V$ of $X$) leads us to the analysis of best simultaneous approximation in Banach spaces. 
To this end
we define the notion of $\mr{F}$-rcp, $\tau$-$\mr{F}$-property-$(P_1)$ and $\tau$-$\mr{F}$-SACP, 
for a triplet $(X,V,\mf{F})$, from \cite{ST1}.
For $N\in\mb{N}$, 
we consider $\mb{R}^N$ endowed with the $\|.\|_\iy$ norm. In this study, $\tau$ denotes the norm or weak topology on $X$.
Denote $\mc{F}(X)$ the set of all nonempty finite subsets of $X$.
For $x\in X$, a $\tau$-closed subset $V$ of $X$, 
$F=\{x_1,\cdots,x_N\}\in \mc{F}(X)$, 
and a function $f:\mb{R}^N_+\ra\mb{R}_{+}$, we define  $r_f(x,F) = f(\|x-x_1\|,\cdots,\|x-x_N\|)$,  $\trr_V^f(F) = \underset{v\in V}{\inf}r_f(v,F)$,
$\tr{Cent}_V^f(F) = \{v\in V:r_f(v,F)=\tr{rad}_V^f(F)\}$ and
 $\de$-$\tr{Cent}_V^f(F) = \{v\in V:r_f(v,F)\leq \tr{rad}_V^f(F)+\de\}$.
\bdfn\label{D1}
Let $V$ be a $\tau$-closed subset of $X$, 
let $\mf{F}\ci\mc{F}(X)$ be a family of finite subsets of $X$, and 
let $\mr{F}\subset\cup_{N\in\mb{N}}\{f:\mb{R}^N_+\to\mb{R}_+\}$ be a collection of functions.
The triplet $(X,V,\mf{F})$ is said to have
\bla
\item  {\it restricted $\mr{F}$-center property} ($\mr{F}$-rcp in short) if 
$\trC_V^f(F)\neq\es$
for all $F=\{x_1,\cdots,x_N\}\in \mf{F}$ 
and 
$f:\mb{R}^N_+\ra\mb{R}_+$ in $\mr{F}$.

\item  {\it $\tau$-$\mr{F}$-property-$(P_1)$} if, $(X,V,\mf{F})$ has $\mr{F}$-rcp and for $F=\{x_1,\cdots,x_N\}\in \mf{F}$, a function $f:\mb{R}^N_+\ra\mb{R}_+$ in $\mr{F}$, and a $\tau$-neighborhood $W$ of $0$, there exists $\de>0$ such that $\de$-$\tr{Cent}_V^f(F)\ci\tr{Cent}_V^f(F)+W$.
\item  {\it $\mr{F}$-simultaneous approximative $\tau$-compactness property} ($\tau$-$\mr{F}$-SACP in short) if for $F=\{x_1,\cdots,x_N\}\in \mf{F}$, a function $f:\mb{R}^N_+\ra\mb{R}_+$ in $\mr{F}$ and a sequence $(v_n)\ci V$ such that $r_f(v_n,F)\ra \textrm{rad}_V^f(F)$ implies that $(v_n)$ has a $\tau$-convergent subsequence.
\el
\edfn
When $\tau$ is the norm topology on $X$, we simply write $\mr{F}$-property-$(P_1)$ (for $(b)$) or $\mr{F}$-SACP (for (c)). We refer to \cite{ST2,ST1} for various examples and characterizations of these properties.
Throughout this paper, we restrict 
the class $\mr{F}$ to the following family of functions 
\[
\mr{F}_{cmc}:=\bigcup_{N\in\mb{N}}\left\{
f:\mb{R}^N_+\ra\mb{R}_+\, \middle|\, 
f \mbox{ is convex, monotone, and coercive}
\right\},
\]
where $f:\mb{R}^N_+\ra\mb{R}_+$ is said to be
\bla
\item {\it monotone} if for $\vp_1,\vp_2\in\mb{R}^N_+$, $\vp_1\leq\vp_2$ implies $f(\vp_1)\leq f(\vp_2)$,
where $\vp_1\leq\vp_2$ if and only if $\vp_1(i)\leq\vp_2(i)$ for all $i=1,\cdots,N$,
\item {\it coercive} if $f(\vp)\ra\iy$ as $\|\vp\|_\iy\ra\iy$.
\el



Let us recall the notions of ideals and $M$-ideals.
\bdfn
\bla
\item \cite{AL} A subspace $Y$ of $X$ is said to be an {\it ideal} of $X$ if $Y^\perp$ is the kernel of a norm-1 projection on $X^*$.
\item \cite{H} A subspace $Y$ of $X$ is said to be an {\it $M$-ideal} in $X$ if $X^*=Y^\perp\bo_{\ell_1}Z$ for some subspace $Z$ of  $X^*$.
\el
\edfn
Several examples and properties of $M$-ideals can be found in \cite{H}. A weakening of $M$-ideals known as {\it semi $M$-ideals} can be found in \cite{AL1}.

Let us recall the notion of {\it polyhedral direct sum} of Banach spaces.

\bdfn\cite{LV}
\bla
\item A function $\pi:\mb{R}^k_+\ra\mb{R}_+$ is said to be a {\it norm} if it is subadditive, positively homogeneous and $\pi (t)=0\Lglra t=0$.
\item A norm $\pi:\mb{R}^k_+\ra\mb{R}_+$ is said to be {\it polyhedral} if $\pi$ is of the form $\pi(t)=\underset{1\leq j\leq l}{\max}g_j(t)$, where $g_1,\cdots,g_l\in(\mb{R}^k)^*$. In this case, we say that $\{g_1,\cdots,g_l\}$ generates $\pi$.
\item $X$ is said to be a {\it polyhedral direct sum} of Banach spaces $X_1,\cdots, X_k$ if $X=X_1\oplus\cdots\oplus X_k$ and the norm on $X$ is of the form $\|x\|_\pi=\pi(\|x(1)\|,\cdots,\|x(k)\|)$, for all $x=(x(1),\cdots,x(k))\in X$, where $\pi$ is a polyhedral nondecreasing norm on $\mb{R}^k_+$. We denote the polyhedral sum of $(X_i)_{i=1}^k$ along with the polyhedral norm $\pi$ by $(\oplus X_i)_\pi $.
\el
\edfn
For Banach spaces $X$ and $Y$,  $\mc{B}(X\times Y)$ denotes the space of all bilinear forms on $X\otm Y$. For $z\in X\otm Y$, let $\sum_{i=1}^nx_i\otm y_i$ be a representation of $z$. We define the {\it cross norm} on $X\otm Y$, by $\la(z)=\sup\{\|\sum_{i=1}^n\phi(x_i)y_i\|:\phi\in S_{X^*}\}$. The completion of $(X\otm Y,\la)$  in $\mc{B}(X\times Y)^*$ is called the {\it injective tensor product} of $X$ and $Y$ and is denoted by $X\hat{\otm}_\e Y$. We refer to \cite{JD} for properties and results related to the tensor product.

Let us recall the notion of K\"othe-Bochner spaces from \cite{PKL}.
\bdfn\label{D3}
Let $(\Om,\Si,\mu)$ be a $\si$-finite measure space. A {\it K\"othe function space} $E$ is a Banach space of real-valued functions satisfying the following.
\bln
\item For any finite measurable set $A\in\Si$, $\chi_A\in E$.
\item For $g\in E$ and a measurable function $f$ such that $|f(\om)|\leq|g(\om)|$ for almost all $\om\in\Om$, implies $f\in E$ and $\|f\|_E\leq\|g\|_E$.
\el
\edfn
A Banach space $E$ is said to be {\it order continuous} if for any decreasing sequence $(f_n)$ in $E$ that converges to 0 pointwise for almost all $\om\in\Om$ implies $\lim_{n}\|f_n\|=0$.

Let $E$ be a Banach space with a 1-unconditional basis $(e_n)$. Let $(X_n)$ be a sequence of Banach spaces. Let us define $(\sum\oplus X_n)_E=\{(x_n):x_n\in X_n~\tr{and}~\sum_{n=1}^\iy\|x_n\|e_n\in E\}$. The norm of an element $(x_n)\in (\sum\oplus X_n)_E$ is defined as: $\|(x_n)\|_{(\sum\oplus X_n)_E}=\big\|\sum_{n=1}^\iy\|x_n\|e_n\big\|_E$. It is easy to check that with this norm $(\sum\oplus X_n)_E$ is a Banach space.

Let $X$ be a Banach space and $E$ be the space as stated in Definition~\ref{D3}. The {\it K\"othe-Bochner function space} $E(X)$ is the collection of all $X$-valued measurable functions $F$ such that $[\om\mapsto \|F(\om)\|_X]\in E$. The norm on $E(X)$ is defined as $\|F\|_{E(X)}=\big\|\|F(.)\|_X\big\|_E$.
One may find a detailed discussion on K\"othe-Bochner function spaces in \cite{PKL}.

\subsection{Overview and Roadmap}

This subsection outlines how our main themes evolve into Sections~2 and 3.
We focus on the transfer of these properties under local constraints. 
Section~2 resolves Rao’s question by showing that, under local constraints, the $M$-ideal property descends from $Z_1\subset X$ to $Z_2=Z_1\cap Y\subset Y$. 
Abstracting from this proof, we formulate a unified transfer mechanism and apply it to obtain stability for central subspaces and \GC, with functorial extensions to polyhedral direct sums, K\"othe-Bochner function spaces, and injective tensor products.

In a complementary direction, Section~3, first, we recall that while sums preserve $M$-ideals (or semi $M$-ideals), centrality and Vesel\'y’s \GC~ may fail on sums; we then develop a framework for best simultaneous approximation based on convex, monotone, coercive functions 
on $\mb{R}^N_+$ (for finite $N$).
Under mild assumptions on $Y$ (e.g., finite dimensionality or reflexivity),
we establish transfer theorems for the restricted center property ($\mr{F}$-rcp), the $\mr{F}$-property-$(P_1)$, and the $\mr{F}$-simultaneous approximative compactness property ($\mr{F}$-SACP) from $Z$ to the closed sum $Y+Z$.


\section{Local-Constraint Transfer --- Main Result and Consequences}

\subsection{Motivation}

In this section, we establish that, under local constraints, 
centrality descends to intersections, and subsequently we settle Rao’s problem for $M$-ideals.

Throughout the section, 
$X$ is a real Banach space, $Y\subset X$ is closed, and
$Z_1\subset X$ is closed; we set $Z_2:=Z_1\cap Y$.
We use the local-constraint hypothesis: for each $z\in Z_1$, there exists
norm-one projections
\[
P:\spn\{z,Z_2\}\to Z_2,\qquad
Q:\spn\{z,Y\}\to Y,
\]
such that $Pz=Qz$.

\subsection{Central subspaces under local constraints}
Our first positive result shows that, under local constraints, 
centrality descends from $Z_1$ to the intersection $Z_2=Z_1\cap Y$.

\begin{thm}\label{T7}
Let $Z_1,Y$ be subspaces of $X$ and $Z_2=Z_1\cap Y$. Suppose $Z_2$ is locally constrained w.r.t. $Z_1$ via $Y$. If $Z_1$ is a central subspace of $X$, then $Z_2$ is a central subspace of $Y$.
\end{thm}
\begin{proof}
Let  $\{z_1,\cdots,z_n\}\ci Z_2$ and $\cap_{i=1}^nB_Y(z_i,r_i)\neq\es$. Then, the balls $(B_X(z_i,r_i))_{i=1}^n$ have a nonempty intersection in $X$ and thus, by our assumption, there exists $z\in Z_1$ such that $z\in\cap_{i=1}^nB_X(z_i,r_i)\cap Z_1$. Since $Z_2$ is locally constrained w.r.t. $Z_1$ via $Y$, there exist two projections $Q:\spn\{z,Y\}\ra Y$ and $P:\spn\{z,Z_2\}\ra Z_2$ such that $Pz=Qz$. Then, for all $i=1,\cdots,n$, we have, $\|Qz-z_i\|=\|P(z-z_i)\|\leq\|z-z_i\|
\leq r_i$. Thus, we have $Qz\in\cap_{i=1}^nB_Y(z_i,r_i)\cap Y$.
\end{proof}

\begin{rem}[One-line transfer]
The argument uses only the following identity
$\|Pz-a\|=\|Q(z-a)\|\le \|z-a\|$ for centers $a\in Y$.
We will exploit this one-line transfer repeatedly below.
\end{rem}

\begin{rem}
Let us observe that under the same assumptions on subspaces $Z_1,Z_2$ and $Y$ as stated in Theorem~\ref{T7}, for a family of subsets $\mc{A}$ (finite, compact, closed and bounded or power set), if  $Z_1$ is an $\mc{A}(Z_1)$-$C$-subspace of $X$, then $Z_2$ is an $\mc{A}(Z_2)$-$C$-subspace of $Y$.
\end{rem}
Let us recall from the introduction that when $Y$ is an almost constrained subspace of $X$, then $Y$ is a locally constrained subspace with respect to $X$ via $Y^{**}$. Using this fact and Theorem~\ref{T7}, we have the following.
\begin{cor}\label{C2}
Let $X$ be a Banach space such that $X\in\GC$ and  $Y$ be an almost constrained subspace of $X$. Then $Y\in\GC$.
\end{cor}

\begin{rem}
Using similar arguments as in Theorem~\ref{C2}, we can show that for a family of subsets $\mc{A}$ (finite, compact, closed and bounded or power set), if $X$ has $\mc{A}(X)$-IP and $Y\ci X$ is an almost constrained subspace of $X$, then Y has $\mc{A}(Y)$-IP.
\end{rem}
The following result characterizes central subspaces of a space in $\GC$ in terms of the geometric structure of the subspace itself.
\begin{thm}
Let $X\in\GC$ and $Y$ be an ideal in $X$. Then $Y$ is a central subspace of $X$ if and only if $Y\in \GC$.
\end{thm}
\begin{proof}
Let $y_1,\cdots,y_n\in Y$ and $r_1,\cdots,r_n>0$ be such that $\cap_{i=1}^nB_{Y^{**}}[y_i,r_i]\neq\es$. Then $\cap_{i=1}^nB_{X^{**}}[y_i,r_i]\neq\es$ and as $Y$ is a central subspace of $X$, we have $\cap_{i=1}^nB_Y[y_i,r_i]\neq\es$. This establishes our assertion.

Conversely, suppose $Y\in\GC$. Let $y_1,\cdots,y_n\in Y$ and $r_1,\cdots,r_n>0$ be such that $\cap_{i=1}^nB_X[y_i,r_i]\neq\es$. Hence $\cap_{i=1}^nB_{X^{**}}[y_i,r_i]\neq\es$. Let $x^{**}\in\cap_{i=1}^nB_{X^{**}}[y_i,r_i]$. As $Y$ is an ideal in $X$, there exists a projection of $P$ on $X^*$ such that $ker(P)=Y^\perp$. Thus $P^*$ is a projection of norm-1 from $X^{**}$ to $Y^{\perp\perp}$. Now $Px^{**}\in\cap_{i=1}^nB_{Y^{\perp\perp}}[y_i,r_i]$. As $Y\in\GC$, we have $\cap_{i=1}^nB_Y[y_i,r_i]\neq\es$. This completes the proof.
\end{proof}
\subsection{$M$-ideals under local constraint (Rao's problem)}
We now pass from centrality to $M$-ideals and answer Rao's question in the affirmative.
It is well-known that the sum of two $M$-ideals is an $M$-ideal (see \cite{H}). In  \cite{Rao}, it was asked whether an $M$-ideal analogue of Lemma~2 from \cite{Rao} is valid under the locally constrained hypothesis. We get an affirmative answer to the problem 
by the same one-line transfer.

\begin{thm}\label{T3}
Let $Y,Z_1$ be subspaces of $X$ and $Z_2=Z_1\cap Y$. Suppose $Z_2$ is locally constrained w.r.t. $Z_1$ via $Y$. Then if $Z_1$ is a (semi) $M$-ideal in $X$, then $Z_2$ is a (semi) $M$-ideal in $Y$.
\end{thm}
\begin{proof}
We prove the result for $M$-ideals. The proof for semi $M$-ideals is similar.

Let $(B_Y(a_i,r_i))_{i=1}^3$ be closed balls in $Y$ such that $\cap_{i=1}^3B_Y(a_i,r_i)\neq\es$ and $B_{Y}(a_i,r_i)\cap {Z_2}\neq\es$ for all $i=1,2,3$. As $Z_1$ is an $M$-ideal in $X$, for $\e>0$, there exists $z_0\in Z_1$ such that $z_0\in\cap_{i=1}^3B_X(a_i,r_i+\e)\cap Z_1$. Let $P,Q$ be projections for $z_0$ as defined in the definition of locally constrained spaces. Then, $Pz_0\in Z_2$ and for all $i=1,2,3$, we have
\[
\|Pz_0-a_i\|=\|Qz_0-a_i\|=\|Q(z_0-a_i)\|\leq\|z_0-a_i\|\leq r_i+\e.
\]
This concludes the proof.
\end{proof}

\subsection{Further Consequences of the Local-Constraint Transfer}
From now on, we will use the same one‑line transfer established above repeatedly. 
We show
polyhedral direct sums, K\"othe-Bochner function spaces, injective tensor
products (including $C(K,\cdot)$), and Vesel\'y's $(\mathrm{GC})$.

\begin{thm}
For $1\leq i\leq k$, let $X_i$ be Banach spaces and let $Y_i$ be subspaces of $X_i$. Let $\pi$ be the corresponding polyhedral norm. Let $Z_{1i}$ be subspaces of $X_i$ and  $Z_{2i}=Y_i\cap Z_{1i}$ and $Z_{2i}$ be locally constrained w.r.t. $Z_{1i}$ via $Y_i$ for all $i=1,\cdots,k$. Further suppose that $X=(X_i)_{\pi}$, $Y=(Y_i)_\pi$, $Z_1=(Z_{1i})_\pi$ and $Z_2=(Z_{2i})_\pi$ are the corresponding polyhedral direct sums of the spaces. If $Z_1$ is a central subspace of $X$, then $Z_2$ is a central subspace of $Y$.
\end{thm}
\begin{proof}
Let $(B_Y(a_i,r_i))_{i=1}^n$ be a collection of $n$ closed balls  with centers in $Z_2$ such that $\cap_{i=1}^nB_Y(a_i,r_i)\neq\es$. Since the balls have a nonempty intersection in $X$, by our hypothesis, the balls have a nonempty intersection in $Z_1$. Let $z_0=(z_0(1),\cdots,z_0(k))\in\cap_{i=1}^nB_X(a_i,r_i)\cap Z_1$. By our hypothesis, there exist norm-1 projections $P_i:\spn\{z_0(i),Z_{2i}\}\ra Z_{2i}$ and $Q_i:\spn\{z_0(i),Y_i\}\ra Y_i$ such that $P_iz_0(i)=Q_iz_0(i)$ for all $i=1,\cdots,k$. Now let us define  $P:\spn\{z_0,Z_2\}\ra Z_2$ by $P(\al z_0+w)=(\al P_1z_0(1)+P_1w(1),\cdots,\al P_kz_0(k)+P_kw(k))$ for all $w\in Z_2$ and $\al\in\mb{R}$ and $Q:\spn\{z_0,Y\}\ra Y$ by $Q(\al z_0+y)=(\al Q_1z_0(1)+Q_1y(1),\cdots,\al Q_kz_0(k)+Q_ky(k))$ for all $y\in Y$ and $\al\in\mb{R}$. Then it is easy to see that $P,Q$ are norm-1 projections and $Pz_0=Qz_0$. Now for all $i=1,\cdots,k$ we have 
\[
\|Pz_0-a_i\|_\pi=\|Qz_0-a_i\|_\pi=\|Q(z_0-a_i)\|_\pi\leq\|z_0-a_i\|_\pi\leq r_i.
\]
Thus $Pz_0\in\cap_{i=1}^nB_Y(a_i,r_i)\cap Z_2$ and this completes the proof.
\end{proof}
We now prove the stability of central subspaces in K\"othe-Bochner function spaces. First, we prove the following.
\begin{thm}
Let $E$ be a Banach space with an 1-unconditional basis $(e_n)$, $(X_n)$ be a sequence of Banach spaces,  and $Y_n$ be subspaces of $X_n$ for all $n\in\mb{N}$. Let $X=(\sum\oplus X_n)_E$ and $Y=(\sum\oplus Y_n)_E$. If $Y_n$ is a central subspace of $X_n$  for all $n\in\mb{N}$, then $Y$ is a central subspace of $X$.
\end{thm}
\begin{proof}
Let $x=(x(n))\in X$ and  $y_1,\cdots,y_k\in Y$. Then, by our assumption, there exists $y_n\in Y_n$ such that $\|y_n-y_i(n)\|_{X_n}\leq\|x(n)-y_i(n)\|_{X_n}$ for all $i=1,\cdots,k$ and $n\in\mb{N}$. Taking 0 as an additional point in $Y$, we have $\|y_n\|_{X_n}\leq\|x(n)\|_{X_n}$ for all $n\in\mb{N}$. Let $y=(y_n)$. Then $y\in Y$ and for all $i=1,\cdots,k$, we have
\[\|y-y_i\|_X=\big\|(\|y_n-y_i(n)\|)_{X_n}\big\|_E\leq \big\|(\|x(n)-y_i(n)\|)_{X_n}\big\|=\|x-y_i\|_X.\] 
This completes the proof.
\end{proof}
\begin{thm}
Let $E$ be an order continuous K\"othe function space and $Y$ be a separable subspace of $X$. If $Y$ is a central subspace of $X$, then $E(Y)$ is a central subspace of $E(X)$.
\end{thm}
\begin{proof}
Let $F\in E(X)$ and $H_1,\cdots,H_k\in E(Y)$. Let $G=\{(t,y)\in\Om\times Y:\|y-H_i(t)\|_X\leq\|F(t)-H_i(t)\|_X~\tr{for all}~i=1,\cdots,k\}=\cap_{i=1}^k\{(t,y)\in(\Om\times Y):\|y-H_i(t)\|_X\leq\|F(t)-H_i(t)\|_X\}$. Since $Y$ is a central subspace of $X$, the projection of $G$ onto the first coordinate is $\Om$.

 Since $F,H_1,\cdots,H_k$ are measurable, $G$ is measurable. Thus, by a consequence of von Neumann selection theorem \cite[Theorem~7.2]{TP}, we get a measurable function $H:\Om\ra Y$ such that $(t,H(t))\in G$ for almost all $t$.
Thus, $\|H(t)-H_i(t)\|_X\leq\|F(t)-H_i(t)\|_X$ for almost all $t$ and so $\|H-H_i\|_{E(X)}\leq\|F-H_i\|_{E(X)}$ for all $i=1,\cdots,k$. This completes the proof.
\end{proof}
We now show that the property of being a central subspace is preserved in K\"othe-Bochner spaces under certain locally constrained assumptions.
\begin{thm}
Let $E$ be an order continuous K\"othe function space.  Let $Y$ be a constrained subspace of $X$ with projection $P$. Let $Z_1,Z_2$ be subspaces of $X$ such that $P(Z_1)\ci Z_1$ and $Z_2=Z_1\cap Y$. If $E(Z_1)$ is a central subspace of $E(X)$, then $E(Z_2)$ is a central subspace of $E(Y)$.
\end{thm}
\begin{proof}
Let $F\in E(Y)$ and $H_1,\cdots,H_k\in E(Z_2)$. By our assumption, there exists $H\in E(Z_1)$ such that $\|H-H_i\|_{E(X)}\leq\|F-H_i\|_{E(X)}$ for all $i=1,\cdots,k$. Let $G=P\circ H$. Then, for all $i=1,\cdots,k$, we have $\|G(t)-H_i(t)\|_X=\|P\circ H(t)-H_i(t)\|_X\leq\|H(t)-H_i(t)\|_X$. Thus $\|G-H_i\|_{E(X)}\leq\|H-H_i\|_{E(X)}\leq\|F-H_i\|_{E(X)}$. Hence, the result follows. 
\end{proof}

\begin{thm}\label{T5}
Let $X,W$ be Banach spaces and $Y$ be a constrained subspace of $X$ with a projection $P$. Suppose $Z_1$ is a subspace of $X$ such that $P(Z_1)\ci Z_1$ and $Z_2=Z_1\cap Y$. Then $W\hat{\otimes}_\e Z_2$ is a central subspace of $W\hat{\otimes}_\e Y$, if $W\hat{\otimes}_\e Z_1$ is a central subspace of $W\hat{\otimes}_\e X$. 
\end{thm}
\begin{proof}
Let $(B(u_i,r_i))_{i=1}^n$ be $n$ closed balls in $W\hat{\otimes}_\e Y$ such that $\cap_{i=1}^nB(u_i,r_i)\cap W\hat{\otm}_\e Y\neq\es$, where $u_i\in W\hat{\otimes}_\e Z_2$ for all $i=1,\cdots,n$. It is easy to see that $I\otm P:W\hat{\otm}_\e X\ra W\hat{\otm}_\e Y$ is a norm-1 projection onto $W\hat{\otm}_\e Y$ and $I\otm P(W\hat{\otm}_\e Z_2)\ci W\hat{\otm}_\e Z_1$. By our assumption, $\cap_{i=1}^n(B(u_i,r_i))\cap W\hat{\otm}_\e Z_1\neq\es$. Let $u\in \cap_{i=1}^n(B(u_i,r_i))\cap W\hat{\otm}_\e Z_1$. Then $I\hat{\otm}P(u)\in W\hat{\otm}_\e Z_2$ and for all $i=1,\cdots,n$, we have $\|I\hat{\otm}P(u)-u_i\|=\|I\hat{\otm}P(u-u_i)\|\leq\|u-u_i\|\leq r_i$. This concludes our assertion.
\end{proof}
As a consequence of Theorem \ref{T5}, we have the following.
\begin{cor}
Let $Y$ be a constrained subspace of $X$  with a projection $P$, and $K$ be a compact Hausdorff space. Suppose $Z_1$ is a subspace of $X$ such that $P(Z_1)\ci Z_1$ and $Z_2=Z_1\cap Y$. Then $C(K,Z_2)$ is a central subspace of $C(K,Y)$,    if $C(K,Z_1)$ is a central subspace of $C(K,X)$.
\end{cor}
\begin{proof}
For a Banach space $W$, using the identification $C(K,W)=C(K)\hat{\otm}_\e W$, the result follows from Theorem~\ref{T5}.
\end{proof}

\section{Simultaneous Approximation on Sums: A Unified Framework}
\subsection{Motivation: Why sums fail}

We first show that generalized centers $\GC$ property can fail to be preserved under sums.
Specifically, we construct subspaces $Y,Z\subset X$ such that $Y,Z\in\GC$ and $Y+Z$ is closed,
yet $Y+Z\notin\GC$. 
\begin{ex}\label{E2}
Let $Y=\{(x_n)\in c_0:\sum_{n=1}^\iy\frac{1}{2^n}x_{2n-1}=\sum_{n=1}^\iy\frac{1}{2^n}x_{2n}=0\}$ and $Z=\spn\{(1,1,0,\cdots)\}$ and $X=Y+Z$. Let $f_1=(0,\frac{1}{2},0,\frac{1}{2^2},\cdots), f_2=(\frac{1}{2},0,\frac{1}{2^2},\cdots)\in\ell_1$. Then $Y=ker(f_1)\cap ker(f_2)$. As $1=\|f_1\|=2|f_1(2)|$ and $1=\|f_2\|=2|f_2(1)|$, it follows from \cite[Theorem~6.3]{BP} that $Y$ is constrained in $c_0$. Hence $Y\in \GC$.  Also, since $Z$ is finite dimensional, $Z\in\GC$. Let us observe that $X=\{(x_n)\in c_0:\sum_{n=1}^\iy\frac{1}{2^n}(x_{2n}-x_{2n-1})=0\}$. Now consider $f=(-\frac{1}{2},\frac{1}{2},-\frac{1}{2^2},\frac{1}{2^2},\cdots)\in \ell_1$. Then $X=ker(f)$, $\|f\|_\iy=\frac{1}{2}$ and $\|f\|_1=2$. Now, as $2\|f\|_\iy<\|f\|_1$ and $supp(f)$ is not finite, it follows  from \cite[Theorem~2]{LV1} that $X\notin \GC$. 
\end{ex}

Central subspaces may also fail to be preserved under sums.
We now give two examples to show that the sum of two central subspaces may not be a central subspace.
\begin{ex}

Let $Y$ and $Z$ be as in Example~\ref{E2}. It is observed that both $Y$ and $Z$ are constrained in $c_0$  and thus they are central subspaces of $c_0$. But as $Y+Z$ is not constrained in $c_0$, it follows from \cite[Proposition~13]{Rao1} that $Y+Z$ is not a central subspace of $c_0$.
\end{ex}
\begin{ex}
Let $X$ be the space $\mb{R}^3$ with $\|.\|_\iy$ norm, $Y_1=\spn\{(1,0,-1)\}$ and $Y_2=\spn\{(0,1,-1)\}$. Then $Y_1+Y_2=\{(x_1,x_2,-x_1-x_2):x_1,x_2\in\mb{R}\}$. Clearly, $Y_1,Y_2$ are constrained in $X$, and hence they are central subspaces of $X$. Now let $x=(-\frac{1}{2},-\frac{1}{2},-\frac{1}{2})$, $y_1=(-2,1,1)$, $y_2=(1,1,-2)$, $y_3=(1,-2,1)$. Now as $\|x-y_i\|=\frac{3}{2}$ for all $i=1,2,3$, we have $\cap_{i=1}^3B_X(y_i,\frac{3}{2})\neq\es$. However, it is easy to see that the balls $(B_X(y_i,\frac{3}{2}))_{i=1}^3$ have no common point of intersection in $Y_1+Y_2$.
\end{ex}
This shows that sums offer no general stability mechanism, leading us to consider them under certain mild assumptions.
\subsection{Framework and standing assumptions}\label{subsec:framework}
We fix a family of finite subsets $\mf F=\mc F(X)$ and a function class
$\mr F =\mr F_{\mathrm{cmc}}$.
In this section, we derive that under certain assumptions, for closed subspaces $Y$ and $Z$ of $X$, such that $(X,Y,\mc{F}(X))$ and $(X,Z,\mc{F}(X))$ have $\mr{F}$-rcp ($\tau$-$\mr{F}$-property-$(P_1)$) ($\tau$-$\mr{F}$-SACP), where $\mr{F}$ is a collection of functions and $Y+Z$ is closed, implies $(X,Y+Z,\mc{F}(X))$ has $\mr{F}$-rcp ($\tau$-$\mr{F}$-property-$(P_1)$) ($\tau$-$\mr{F}$-SACP).

We state two results that are essential for deriving several consequences in this section.
\begin{thm}\cite[Theorem~6.20]{WR}\label{T8}
Let $Y$ and $Z$ be two subspaces of $X$ such that $X=Y+Z$. Then there exists $\ga>0$ such that every $x\in X$  has a representation $x=y+z$, and $\|y\|+\|z\|\leq\ga\|x\|$.
\end{thm}
\begin{prop}\cite[Proposition~1.6]{RRP}\label{T9}
Let $D$ be an open convex subset of $X$ and $f:D\ra\mb{R}$ be a convex function. If $f$ is continuous at $x_0\in D$, then $f$ is locally Lipschitzian at $x_0$, i.e., there exists $M>0$ and $\de>0$ such that $B(x_0,\de)\ci D$ and $|f(x)-f(y)|\leq M\|x-y\|$ for all $x,y\in B(x_0,\de)$.
\end{prop}

\subsection{Main transfer on sums: the $\mr F_{\mathrm{cmc}}$-rcp case}
We first establish the transfer of the restricted $\mr F_{\mathrm{cmc}}$-center property from $Z$ to the closed sum $Y+Z$.

\begin{thm}\label{T6}
Let $Y$ be a reflexive subspace of $X$ and $Z$ be a subspace of $X$ such that $(X,Z,\mc{F}(X))$ has $\mr{F}_{cmc}$-$\ercp$. Then $(X,Y+Z,\mc{F}(X))$ has $\mr{F}_{cmc}$-$\ercp$.
\end{thm}
\begin{proof}
Let $F=\{x_1,\cdots,x_N\}\in\mc{F}(X)$, $f:\mb{R}^N_+\ra\mb{R}_+$ be a function in $\mr{F}_{cmc}$. Let $(y_n)$ be a sequence in $Y$ and $(z_n)$ be a sequence in $Z$ such that $r_f(y_n+z_n,F)\ra\tr{rad}_{Y+Z}^f(F)$. Without loss of generality, we may assume that $(r_f(y_n+z_n,F))$ is a decreasing sequence converging to $\trr_{Y+Z}^f(F)$. Since $f$ is coercive, $(y_n+z_n)$ is bounded, and from Theorem~\ref{T8}, both $(y_n)$ and $(z_n)$ are bounded. Since $Y$ is reflexive, $(y_n)$ has a weakly convergent subsequence, also denoted by $(y_n)$, converging to $y_0\in Y$. Now there exists a sequence of convex combinations $\Ti{y}_n=\sum_{i\in I_n}\la_iy_i$, where $I_n=\{i:p_n<i\leq p_{n+1}\}$, $(p_n)$ is an increasing sequence of natural numbers, $\la_i>0$ for all $i\in I_n$ and $\sum_{i\in I_n}\la_i=1$, such that $\|\Ti{y}_n-y_0\|\ra0$. Let us denote $\Ti{z}_n=\sum_{i\in I_n}\la_iz_i$ for all $n\in\mb{N}$. Now 
\beqa
\trr_{Y+Z}^f(F)&\leq& f(\|\Ti{y}_n+\Ti{z}_n-x_1\|,\cdots,\|\Ti{y}_n+\Ti{z}_n-x_N\|)\\
&=&f(\|\sum_{i\in I_n}\la_i(y_i+z_i-x_1)\|,\cdots,\|\sum_{i\in I_n}\la_i(y_i+z_i-x_N)\|)\\
&\leq& f(\sum_{i\in I_n}\la_i\|y_i-z_i-x_1\|,\cdots,\sum_{i\in I_n}\la_i\|y_i-z_i-x_N\|)\\
&\leq&\sum_{i\in I_n}\la_if(\|y_i+z_i-x_1\|,\cdots,\|y_i+z_i-x_N\|)\\
&\leq& f(\|y_{p_n+1}+z_{p_n+1}-x_1\|,\cdots,\|y_{p_n+1}+z_{p_n+1}-x_N\|)\\
&\ra&\trr_{Y+Z}^f(F).
\eeqa
Thus $r_f(\Ti{y_n}+\Ti{z_n},F)\ra\trr_{Y+Z}^f(F)$. Now by the continuity of $f$ on bounded subsets of $\mb{R}^N_+$, we have 
\beqa
r_f(y_0+\Ti{z_n},F)&=&f(\|y_0+\Ti{z_n}-x_1\|,\cdots,\|y_0+\Ti{z_n}-x_N\|)\\
&\leq& f(\|\Ti{y_n}+\Ti{z_n}-x_1\|+\|\Ti{y_n}-y_0\|,\cdots,\|\Ti{y_n}+\Ti{x_n}-x_N\|\\
&&+\|\Ti{y_n}-y_0\|)\\
&\ra&\trr_{Y+Z}^f(F).
\eeqa 
Let $z_0\in Z$ be such that $r_f(z_0,F-\{y_0\})=\trr_Z^f(F-\{y_0\})$. Now $r_f(y_0+z_0,F)=r_f(z_0,F-\{y_0\})\leq r_f(\Ti{z_n},F-\{y_0\})= r_f(y_0+\Ti{z_n},F)\ra\trr_{Y+Z}^f(F)$. This concludes our assertion.
\end{proof}

As a consequence of Theorem~\ref{T6}, we have the following.
\begin{cor}\label{C1}
Let $Y$ be a reflexive subspace of $X$ and $Z$ be a subspace of $X$ such that $Y+Z$ is closed and $(X,Z,\mc{F}(X))$ has $\mr{F}_{cmc}$-$\ercp$. Then $Y+Z\in\GC$.
\end{cor}

\begin{rem}[Common skeleton on sums]\label{rem:common-skeleton}
All proofs below follow the same pattern as Theorem~\ref{T6}:
\begin{description}
    \item[Step 1] Coercivity
    \item[Step 2] Decomposition with controlled size by Theorem~\ref{T8}
  \item[Step 3] Compactness on $Y$ (finite-dim. or reflexive)
 \item[Step 4]
 Stability from convexity and coordinatewise monotonicity
  \item[Step 5]
Push the $Z$-side hypothesis forward to the sum
\end{description}
\end{rem}

\begin{thm}\label{T2}
Let $Y$ be a finite dimensional subspace of $X$ and $Z$ be a subspace of $X$ such that $(X,Z,\mc{F}(X))$ has $\tau$-$\mr{F}_{cmc}$-property-$(P_1)$. Then $(X,Y+Z,\mc{F}(X))$ has $\tau$-$\mr{F}_{cmc}$-property-$(P_1)$.
\end{thm}
\begin{proof}
Let $F=\{x_1,\cdots,x_N\}\in\mc{F}(X)$ and $f:\mb{R}^N_+\ra\mb{R}_+$ be a function in $\mr{F}_{cmc}$. Let $(y_n)$ and $(z_n)$ be sequences in $Y$ and $Z$ respectively, such that $r_f(y_n+z_n,F)\ra\trr_{Y+Z}^f(F)$. Since $f$ is coercive $(y_n+z_n)$ is bounded. Now from Theorem~\ref{T8}, $(y_n)$ and $(z_n)$ is bounded.  Since $Y$ is finite dimensional, $(y_n)$ has a convergent subsequence $(y_{n_k})$ converging to $y_0\in Y$. Now, as $\|z_{n_k}+y_0-x_i\|\leq\|z_{n_k}+y_{n_k}-x_i\|+\|y_{n_k}-y_0\|$ for all $i=1,\cdots,N$, by the monotonicity of $f$, we have  for all $k\in\mb{N}$,
\begin{align}\label{e1}
f(\|z_{n_k}&+y_0-x_1\|,\cdots,\|z_{n_k}+y_0-x_N\|) \nonumber \\
&\leq f(\|z_{n_k}+y_{n_k}-x_1\|+\|y_{n_k}-y_0\|,\cdots,\nonumber \\
&\qquad \|z_{n_k}+y_{n_k}-x_N\|+\|y_{n_k}-y_0\|)
\end{align}

\beqa
\mbox{Now~}\|(\|z_{n_k}+y_{n_k}-x_1\|+\|y_{n_k}-y_0\|,\cdots,\|z_{n_k}+y_{n_k}-x_N\|+\|y_{n_k}-y_0\|)\\-(\|z_{n_k}+y_{n_k}-x_1\|,\cdots,
\|z_{n_k}+y_{n_k}-x_N\|)\|_\iy=\|y_{n_k}-y_0\|\ra0.
\eeqa 
Since from Proposition~\ref{T9}, $f$ is Lipschitz continuous on bounded subsets of $\mb{R}_+^N$. Hence, for a given $\de>0$ we have,
\begin{align*}
|f(\|z_{n_k}+y_{n_k}-x_1\|+\|y_{n_k}-y_0\|,\cdots,\|z_{n_k}+y_{n_k}-x_N\|+\|y_{n_k}-y_0\|)\\-\trr_{Y+Z}^f(F)|\\
\leq|f(\|z_{n_k}+y_{n_k}-x_1\|,\cdots,\|z_{n_k}+y_{n_k}-x_N\|)-\trr_{Y+Z}^f(F)|+\de,
\end{align*}
for sufficiently large $k$. As $r_f(y_{n_k}+z_{n_k},F)\ra \trr_{Y+Z}^f(F)$, we get,
\begin{equation*}
    \begin{aligned}
      &|f(\|z_{n_k} + y_{n_k} - x_1\| + \|y_{n_k} - y_0\|, \dots, \|z_{n_k} + y_{n_k} - x_N\| + \|y_{n_k} - y_0\|) \\
      &\quad - \trr_{Y+Z}^f(F)| \ra 0.
    \end{aligned}
  \end{equation*}

Thus, from (\ref{e1}) we have 
\[f(\|z_{n_k}+y_0-x_1\|,\cdots,\|z_{n_k}+y_0-x_N\|)\ra\trr_{Y+Z}^f(F).\]
\beqa\mbox{Hence,}~
\trr_Z^f(F-\{y_0\})
&\leq &\lim_k f(\|z_{n_k}+y_0-x_1\|,\cdots,\|z_{n_k}+y_0-x_N\|)\\
&=&\trr_{Y+Z}^f(F)\leq\trr_Z^f(F-\{y_0\}).
\eeqa

Thus, $r_f(z_{n_k},F-\{y_0\})\ra\trr_Z^f(F-\{y_0\})$ and $\trr_{Y+Z}^f(F)=\trr_Z^f(F-\{y_0\})$. Let $W$ be a $\tau$-neighbourhood of 0 and $W_1,W_2$ be two $\tau$-neighbourhoods of 0 such that $W_1+W_2\ci W$. Now, by our assumption, there exists $\de>0$ such that $\de$-$\trC_Z^f(F-\{y_0\})\ci\trC_Z^f(F-\{y_0\})+W_1$. Now there exists $N_1\in\mb{N}$ such that $z_{n_k}\in\de-\trC_{Y+Z}^f(F)$ for all $k\geq N_1$. Thus there exists a sequence $(s_k)_{k\geq N_1}\ci\trC_Z^f(F-\{y_0\})$ such that $z_{n_k}-s_k\in W_1$ for all $k\geq N_1$. Now as $r_f(s_k,F-\{y_0\})=\trr_Z^f(F-\{y_0\})=\trr_{Y+Z}^f(F)$, we have the sequence $(s_k+y_0)\ci\trC_{Y+Z}^f(F)$. Since $y_{n_k}\ra y_0$, there exists $N_2\in\mb{N}$ such that $y_{n_k}-y_0\in W_2$ for all $k\geq N_2$. Thus, $y_{n_k}+z_{n_k}-s_k-y_0\in W_1+W_2\ci W$ for all $k\geq \max\{N_1,N_2\}$.  Hence, $y_{n_k}+z_{n_k}\in\trC_{Y+Z}^f(F)+W$ for all $k\geq\max\{N_1,N_2\}$ and this concludes our assertion.
\end{proof}
\begin{thm}\label{T1}
Let $Y$ be a finite dimensional subspace of $X$, and $Z$ be a subspace of $X$ such that $(X,Z,\mc{F}(X))$ has $\tau$-$\mr{F}_{cmc}$-$\eSACP$. Then, $(X,Y+Z,\mc{F}(X))$ has $\tau$-$\mr{F}_{cmc}$-$\eSACP$. 
\end{thm}

\begin{proof}
Let $F=\{x_1,\cdots,x_N\}\in\mc{F}(X)$, $f:\mb{R}^N_+\ra\mb{R}_+$ be a function in $\mr{F}_{cmc}$. Suppose $(y_n)\ci Y$, $(z_n)\ci Z$ be such that $r_f(y_n+z_n,F)\ra\trr_{Y+Z}^f(F)$. Since $f$ is coercive, $(y_n+z_n)$ is bounded, and hence from Theorem~\ref{T8}, $(y_n)$ and $(z_n)$ are bounded. Since $Y$ is finite dimensional, $(y_n)$ has a convergent subsequence $(y_{n_k})$ converging to $y_0\in Y$. Now proceeding in a similar way as in the proof of Theorem \ref{T2}, we have $r_f(z_{n_k},F-\{y_0\})\ra\trr_Z^f(F-\{y_0\})$ and so, by our assumption, $(z_{n_k})$ has a $\tau$-convergent subsequence $(z_{n_l})$. Thus $(y_n+z_n)$ has a $\tau$-convergent subsequence $(y_{n_l}+z_{n_l})$, and this establishes our assertion. 
\end{proof}

The following example is motivated by the construction in \cite[Theorem~3]{IP}.
\begin{ex}
Consider the canonical basis $(e_n)\ci\ell_1$. Let $Y_n=\spn\{e_1,\cdots,e_n\}$ for all $n\in\mb{N}$. Then $d(e_{n+1},Y_n)=1$ for all $n\in\mb{N}$. Let $U=\spn\{e_1\}$ and $V=\cup_{n=2}^\iy\{(n-1)e_1+e_n\}$. Then $U+V=\cup_{n=2}^\iy\{e_n\}+\spn\{e_1\}$. It is easy to see that $U+V$ is $\tau$-closed. Since $U$ is finite dimensional, $(X,U,\mc{F}(X))$ has $\tau$-$\mr{F}_{cmc}$-$\eSACP$. Also, since all possible  norm-convergent sequences in $V$ are eventually constant, $(X,V,\mc{F}(X))$ has $\tau$-$\mr{F}_{cmc}$-$\eSACP$. Now, as $d(e_{n+1},Y_n)=1$ for all $n\in\mb{N}$, $d(0,U+V)\geq d(0,e_{n+1}+\spn\{e_1\})\geq 1$. Thus we have $\|e_n\|=1=d(0,U+V)$. But the sequence $(e_n)$ has no $\tau$-convergent subsequence in $\ell_1$. Thus $(X,U+V,\mc{F}(X))$ does not have $\tau$-$\mr{F}_{cmc}$-$\eSACP$.
\end{ex}
\begin{rem}
We do not know whether there exist two subspaces $Y$ and $Z$ of $X$ such that $(X,Y,\mc{F}(X))$ and $(X,Z,\mc{F}(X))$ has $\tau$-$\mr{F}_{cmc}$-$\eSACP$, $Y+Z$ is closed but $(X,Y+Z,\mc{F}(X))$ does not have $\tau$-$\mr{F}_{cmc}$-$\eSACP$.
\end{rem}
Let $Y$ be a subspace of $X$. Let $Z$ be a subspace of $Y$ that is finite codimensional in $Y$. Suppose that $\{y_1+Z,\cdots,y_n+Z\}$ is a basis for the quotient space $Y/Z$, where $\{y_1,\cdots,y_n\}\ci Y$. Then, $Y=Z+\spn\{y_1,\cdots,y_n\}$.
\begin{thm}
Let $Z$ be a finite codimensional subspace of $X$. If $(X,Z,\mc{F}(X))$ has $\tau$-$\mr{F}_{cmc}$-property-$(P_1)$ ($\mr{F}_{cmc}$-$\ercp$) ($\tau$-$\mr{F}_{cmc}$-$\eSACP$), then, for all subspaces $Y$ of $X$ with $Z\ci Y\ci X$, $(X,Y,\mc{F}(X))$ has $\tau$-$\mr{F}_{cmc}$-property-$(P_1)$ ($\mr{F}_{cmc}$-$\ercp$) ($\tau$-$\mr{F}_{cmc}$-$\eSACP$).
\end{thm}
\begin{proof}
We give the proof when $(X,V,\mc{F}(X))$ has $\tau$-$\mr{F}_{cmc}$-property-$(P_1)$, as the proof for $\tau$-$\mr{F}_{cmc}$-SACP is similar.

Suppose that $Y$ is a subspace of $X$ such that $Z\ci Y\ci X$. Since $Z$ is finite codimensional in $Y$, there exists a finite dimensional subspace $S$ of $X$ such that $Y=Z+S$. As $S$ is finite dimensional and $(X,Z,\mc{F}(X))$ has $\tau$-$\mr{F}_{cmc}$-property-$(P_1)$, it follows from Theorem \ref{T1} that $(X,Y,\mc{F}(X))$ has $\tau$-$\mr{F}_{cmc}$-property-$(P_1)$.
\end{proof}

\begin{thm}
 Let $Y$ be a reflexive subspace of $X$ and $Z$ be a subspace of $X$ such that $Y+Z$ is closed and $(X,Z,\mc{F}(X))$ has weak-$\mr{F}_{cmc}$-$\eSACP$. Then, $(X,Y+Z,\mc{F}(X))$ has weak-$\mr{F}_{cmc}$-$\eSACP$.   
\end{thm}
\begin{proof}
Let $F=\{x_1,\cdots,x_N\}\in\mc{F}(X)$, $f:\mb{R}^N_+\ra\mb{R}_+$ be a function in $\mr{F}_{cmc}$. Let $(y_n)$ be a sequence in $Y$ and $(z_n)$ be a sequence in $Z$ such that $r_f(y_n+z_n,F)\ra\tr{rad}_{Y+Z}^f(F)$. Without loss of generality, we may assume that $(r_f(y_n+z_n,F))$ is a decreasing sequence converging to $\trr_{Y+Z}^f(F)$. Since $f$ is coercive, $(y_n+z_n)$ is bounded, and from Theorem~\ref{T8}, both $(y_n)$ and $(z_n)$ are bounded. Since $X$ is reflexive, $(y_n)$ has a weakly convergent subsequence, also denoted by $(y_n)$, converging to $y_0\in Y$. We only need to show that $(z_n)$ has a weakly convergent subsequence. 

Suppose $(z_n)$ has no weakly convergent subsequence. Then by \cite[Theorem~1]{RJ}, there exists a $\theta>0$ and a subsequence of $(z_n)$, also denoted by $(z_n)$ such that $d(co(z_1,\cdots,z_k),co(z_{k+1},z_{k+2},\cdots))\geq\theta$ for all $k\in\mb{N}$. Now there exists a sequence of convex combinations $\Ti{y}_n=\sum_{i\in I_n}\la_iy_i$, where $I_n=\{i:p_n<i\leq p_{n+1}\}$, $(p_n)$ is an increasing sequence of natural numbers, $\la_i>0$ for all $i\in I_n$ and $\sum_{i\in I_n}\la_i=1$, such that $\|\Ti{y}_n-y_0\|\ra0$. Let us define $\Ti{z}_n=\sum_{i\in I_n}\la_iz_i$ for all $n\in\mb{N}$. 

Now, as done in the proof of Theorem~\ref{T6}, we have $r_f(\Ti z_n+\Ti y_n,F)\ra\trr_{Y+Z}^f(F)$.
Thus, using the continuity of $f$ on bounded subsets of $\mb{R}^N_+$, we have $r_f(\Ti{z}_n,F-\{y_0\})=f(\|\Ti{z}_n+y_0-x_1\|,\cdots,\|\Ti{z}_n+y_0+x_N\|)\leq f(\|\Ti{z}_n+\Ti{y}_n-x_1\|+\|\Ti{y}_n-y_0\|,\cdots,\|\Ti{z}_n+\Ti{y}_n-x_N\|+\|\Ti{y}_n-y_0\|)\ra\trr_{Y+Z}^f(F)\leq\trr_Z^f(F-\{y_0\})$. By our assumption, $(\Ti{z}_n)$ has a weakly convergent subsequence $(\Ti{z}_{n_k})$. Now $\Ti{z}_{n_k}\in co(z_{{p_{n_k}}+1},\cdots,z_{p_{{n_k}+1}})\ci co(z_1,z_2,\cdots,z_{p_{n_k+1}})$ for all $k\in\mb{N}$. Thus, $d(co(\Ti{z}_{n_1},\cdots,\Ti{z}_{n_k}),co(\Ti{z}_{n_{k+1}},\Ti{z}_{{n_{k+2}}},\cdots))\geq d(co(z_1,\cdots,z_{p_{n_k+1}}),co(z_{p_{{n_k+1}}+1},\cdots))\geq\theta$ for all $k\in\mb{N}$, a contradiction to the fact that $(\Ti{z}_{n_k})$ is weakly convergent. This contradiction establishes our assertion.
\end{proof}

\section*{Declaration}
The author of this work declare that they have no conflict of interest.

\Addresses
\end{document}